\documentclass
[12pt,a4paper]{article}
\usepackage[cp1251]{inputenc}
\usepackage{amssymb, amsfonts, amsmath, latexsym}
\usepackage[english]{babel}

\begin{document}

\Large
\centerline{\bf  On asymorphisms of finitary coarse spaces  }\vspace{6 mm}

\normalsize\centerline{\bf  I.V. Protasov }\vspace{6 mm}


{\bf Abstract.} We characterize finitary coarse spaces $X$ such that every permutation of $X$ is an asymorphism.

\vspace{6 mm}


\vspace{3 mm}



\section{Introduction and results}

Given a set $X$, a family $\mathcal{E}$ of subsets of $X\times X$ is called a {\it coarse structure }  on $X$  if
\vskip 7pt

\begin{itemize}
\item{}   each $E\in \mathcal{E}$  contains the diagonal  $\bigtriangleup _{X}$,
$\bigtriangleup _{X}= \{(x,x)\in X: x\in X\}$;
\vskip 5pt

\item{}  if  $E$, $E^{\prime} \in \mathcal{E}$ then $E\circ E^{\prime}\in\mathcal{E}$ and
$E^{-1}\in \mathcal{E}$,   where    $E\circ E^{\prime}=\{(x,y): \exists z((x,z) \in  E,  \   \ (z, y)\in E^{\prime})\}$,   $E^{-1}=\{(y,x): (x,y)\in E\}$;
\vskip 5pt

\item{} if $E\in\mathcal{E}$ and $\bigtriangleup_{X}\subseteq E^{\prime}\subseteq E  $   then
$E^{\prime}\in \mathcal{E}$;
\vskip 5pt

\item{}   $\bigcup  \mathcal{E}= X\times X $.

\end{itemize}
\vskip 7pt

A subfamily $\mathcal{E}^{\prime} \subseteq \mathcal{E}$  is called a
{\it base} for $\mathcal{E}$  if,
 for every $E\in \mathcal{E}$, there exists
  $E^{\prime}\in \mathcal{E}^{\prime}$  such  that
  $E\subseteq E ^{\prime}$.
For $x\in X$,  $A\subseteq  X$  and
$E\in \mathcal{E}$, we denote
$$E[x]= \{y\in X: (x,y) \in E\},
 \   E [A] = \bigcup_{a\in A}   \    
  E[a]$$
 and say that  $E[x]$
  and $E[A]$
   are {\it balls of radius $E$
   around} $x$  and $A$.

\vskip 10pt
The pair $(X,\mathcal{E})$ is called a {\it coarse space}  
 \cite{b6} 
 or a ballean
\cite{b4}, 
\cite{b5}.

For a coarse   space $(X,\mathcal{E})$, a  subset $B \subseteq X$   is called {\it bounded} if $B \subseteq  E[x]$  for some
$E\in \mathcal{E}$ and $x\in X$.
The family $\mathcal{B}_{(X, \mathcal{E})}$
 of all bounded subsets of  $(X, \mathcal{E})$
is called the {\it bornology} of $(X, \mathcal{E})$.
We recall that a family $ \mathcal{B} $ of subsets of a set $X$ is a bornology if $ \mathcal{B} $ is closed under taking subsets and finite unions, and $ \mathcal{B} $ contains all finite subsets of $X$. 

\vskip 7pt

Let $(X, \mathcal{E})$, $(X^\prime, \mathcal{E}^\prime)$ be coarse spaces. A mapping $f: X\longrightarrow X^\prime$ is called

\vskip 7pt

\begin{itemize}
\item{} {\it bornologous} if $f(B)\in \mathcal{B}_{(X^\prime, \mathcal{E}^\prime)}$ for each $B\in \mathcal{B}_ {(X, \mathcal{E})}$;

\vskip 7pt

\item{} {\it macro-uniform} if, for each $E\in \mathcal{E}$, there exists $E^\prime \in \mathcal{E}^\prime$ such that, for all $x,y\in X$, $(x,y)\in E$ implies $(f(x), f(y)) \in E^\prime$;

\vskip 7pt

\item{} {\it asymorphism} if $f$ s a bijection and $f, f^{-1}$ are macro-uniform.

\vskip 7pt
\end{itemize}

We recall that a coarse space $(X, \mathcal{E})$ is 
{\it discrete} (or {\it thin}) if, for each $E\in \mathcal{E}$, there exists $B\in \mathcal{B}_{(X,\mathcal{E})}$ such that $E[x]= \{ x \}$
for each $x\in X \setminus B$.
Every bornology $\mathcal{B}$ on a set $X$ defines the discrete coarse space $(X, \mathcal{E}_\mathcal{B})$
with the base $\{ E_B : B\in \mathcal{B} \}$, were 
$E_B [x] = B$ if $x\in B$, and  $E_B [x] = \{ x \}$ if 
$x\in X\setminus B$. 
Every discrete coarse space $(X, \mathcal{E})$ coincides with  
$(X, \mathcal{E}_\mathcal{B} )$ for $\mathcal{B}= \mathcal{B}_{(X, \mathcal{E})}$.

\vskip 5pt

For different characterizations of discrete coarse spaces, see Theorem 2.2 in \cite{b1}.

\vskip 10pt

{\bf Theorem 1.} {\it Every bornologous mapping 
$f: X\longrightarrow X$
of a coarse space $(X, \mathcal{E})$ is macro-uniform if and only if $(X, \mathcal{E})$ is discrete.
\vskip 7pt
}
A coarse space 
$(X, \mathcal{E})$ is called 

 \vskip 7pt

\begin{itemize}
\item{} {\it locally finite} if each ball $E[x]$ is finite, equivalently, 
$\mathcal{B}_{(X, \mathcal{E})} =
[X]^{<\omega} $; 

\vskip 7pt

\item{} {\it finitary } if, for each 
$E\in \mathcal{E}$ there exists a natural number $n$ such that $|E [x]|< n$ for each $x\in X$.

\end{itemize}

\vskip 10pt

Let $G$ be a transitive group of permutations of a set $X$. We denote by $X_G$ the set $X$ endowed with the coarse structure with the base
$$ \{\{(x, gx): g\in F \} : F \in [G]^{< \omega},  \  \ id\in F \}.  $$
By [2, Theorem 1], for every finitary coarse structure 
$(X, \mathcal{E})$, there exists a transitive group $G$
 of permutations  of $X$  such that $(X, \mathcal{E})=X_G$.
 For more general results, see \cite{b3}. 
 
\vskip 7pt 
Let $X$ be a set, $\kappa$ be a cardinal, 
$\mathcal{S}_X$ denotes the group of all permutations of $X$,  
$\mathcal{S}_X ^{<\kappa} = \{ g\in \mathcal{S}_X : |supp \ g| < \kappa \}$, 
$supp \ g = \{ x\in X : g(x) \neq x \}$.

\vskip 10pt

{\bf Theorem 2.} {\it Let $(X, \mathcal{E})$ be an infinite finitary coarse space. Then the following statements are equivalent 

\vskip 7pt
$(i)$  every permutation of $X$ is an asymorphism of 
$(X, \mathcal{E})$;

\vskip 7pt
$(ii)$  there exists an infinite cardinal $\kappa$,
$\kappa \leq |X|^+$ such that $(X, \mathcal{E})= X_G$
for $G=\mathcal{S}_X ^{< \kappa}$.}

\vskip 10pt

{\bf Open problem.} {\it Characterize locally finite coarse spaces $(X, \mathcal{E})$ such that every permutation of $X$ is an asymorphism.}

\section{Proofs}

{\it Proof of Theorem 1.} 
\vskip 10pt

Let $(X, \mathcal{E})$
be a discrete coarse space defined by a bornology $\mathcal{B}$  and let $f: X\longrightarrow X$
is bornologous.
We take an arbitrary $B\in \mathcal{B}$ and note that 
$f(E_B [x])\subseteq E_{f(B)} [f(x)]$
for each $x\in X$, so $f$ is macro-uniform.

On the other hand, let $(X, \mathcal{E})$ is not discrete. 
Then there exists $E\in \mathcal{E}$ such that, for each bounded subset $B$ of $(X, \mathcal{E})$, 
one can find $x\in X\setminus B$ such that 
$|E[x] | > 1$.
 Therefore, for some ordinal $\lambda$, we can choose 
 inductively two injective $\lambda$-sequences 
 $(x_\alpha)_{\alpha<\lambda}$,
  $(y_\alpha)_{\alpha<\lambda}$
such that the set $\{ x_\alpha : \alpha< \lambda \}$
is unbounded, 
$y_\alpha)\in E [x_\alpha]$, $\alpha< \lambda$ and 
$y_\alpha \neq x_\beta$ for all $\alpha, \beta < \lambda$.

We define a mapping $f: X\longrightarrow X$ by $f(y_\alpha)= y_0$ for each $\alpha<\lambda$, and 
$f(x)=x$ for each $x\in X\setminus \{ y_\alpha : \alpha<\lambda\}$.
Clearly, $f$ is bornologous. Since $\{ x_\alpha : \alpha<\lambda\}$ is unbounded, $f(x_\alpha)=x_\alpha$,
$y_\alpha \in E [x_\alpha]$ and 
$f(y_\alpha)=y_0 $ for each $\alpha<\lambda$, we conclude that $f$ is not macro-uniform. 
$ \ \  \  \Box $

\vskip 7pt

{\it Proof of Theorem 2.}
\vskip 10pt

$(i)\Longrightarrow (ii)$. We say that a permutation $g$ of $X$ is {\it compatible} with $\mathcal{E}$ if there exists $E\in \mathcal{E}$ such that $(x, gx)\in E$ for each $x\in X$. We note that the set  $G$ of all permutations compatible with $\mathcal{E}$ is a subgroup of $\mathcal{S}_X$ and, by Theorem 1 from \cite{b2},  $(X, \mathcal{E})= X_G$. 

We say that a subset $Y$ of $X$ is crowded if there exists $E\in \mathcal{E}$ such  that $|E[y]|>1 $ for each $y\in Y$. We take the minimal cardinal $\kappa$,
$\kappa\leq |X|^+$ such that, for each $\lambda<\kappa$,
$(X, \mathcal{E})$ has a crowded subset of cardinality 
$\lambda$.

We show that $G=\mathcal{S} _X ^{<\kappa}$.
If $g\in G$ then 
$|supp \ g|<\kappa$ 
because the set $supp \ g$ is crowded, so $g\in \mathcal{S} _X ^{<\kappa}$ and $G\subseteq \mathcal{S} _X ^{<\kappa}$.
\vskip 7pt

To prove $\mathcal{S} _X ^{<\kappa}\subseteq G$, we need the following auxiliary statement. 
\vskip 7pt

$(\ast)$  {\it Let $Y$ be a subset of $X$ such that 
$|Y|=|X|$, $\lambda$ be a cardinal, $\lambda<\kappa$.
Then there exists two injective $\lambda$-sequences 
$(x_\alpha) _ {\alpha<\lambda}$, 
$(y_\alpha) _ {\alpha<\lambda}$ 
in $Y$
and $H\in \mathcal{E}$
such that $\{x_\alpha : \alpha<\lambda \} \cap \{ y_\alpha : \alpha<\lambda \}=\emptyset $
and $(x_\alpha, y_\alpha)\in H$ for each $\alpha < \lambda$. }

\vskip 10pt

By the choice of $\kappa$, we can choose $E\in  \mathcal{E}$ and injective $\lambda$-sequences 
$(a_\alpha ) _ {\alpha<\lambda}$, 
$(b_\alpha ) _ {\alpha<\lambda}$
such that 
$\{ a_\alpha : \alpha<\lambda\} \  \bigcap
 \ \
\{b_\alpha :  \alpha<\lambda \}=\emptyset$
and
$(a_\alpha, b_\alpha ) \in E$ for each 
$\alpha<\lambda$.
Passing to subsequences, we may suppose that $|Y\setminus \{ a_\alpha, \ b_\alpha:  \alpha<\lambda \} |= |X|.$
We choose two injective $\lambda$-sequences 
$(x_\alpha)_{ \alpha< \lambda}$, 
$(y_\alpha)_{\alpha< \lambda}$
in $Y\setminus \{ a_\alpha, \ b_\alpha:  \alpha<\lambda \}$
such that 
$\{ x_\alpha: \alpha< \lambda \} \ \cap$
$\{ y_\alpha: \alpha< \lambda \}= \emptyset$.
Then we define an involution $f$ of $X$ by 
$f a_\alpha =x_\alpha$, $f b_\alpha =y_\alpha$ and $fx = x$ for each $x\in X\setminus \{a_\alpha, \ b_\alpha, \  x_\alpha, \ y_\alpha : \alpha<\lambda \}$.
Since $f$ is macro-uniform, there exists $H\in \mathcal{E}$ such that $(x,y)\in E$ implies $(fx, fy)\in H$. Hence, $(x_\alpha, y_\alpha)\in H$ for each $\alpha<\lambda$.

\vskip 10pt

Now let $g\in \mathcal{S}_X^{<\kappa}$,
$A=supp \ g$. We prove that $g$ is compatible with $\mathcal{E}$, so $g\in G$. By the 3-Sets Lemma, there exists a partition $A_1 , A_2 , A_3$  of $A$ such that 
$A_i\bigcap gA_i = \emptyset$, $i\in \{1, 2, 3 \}$.
We suppose that $|X\setminus (A_1 \bigcup g A_1))| = |X|$, denote $Y=|X\setminus (A_1 \bigcap gA_1)|$,
$\lambda=|A_1|$ and apply $(\ast)$ to choose corresponding
$(x_\alpha)_{\alpha<\lambda}$, $(y_\alpha)_{\alpha<\lambda}$
and $H\in \mathcal{E}$. We enumerate $A_1 = \{ a_\alpha : \alpha<\lambda \}$ and define an involution $h$ of $X$ by $hx_\alpha = a_\alpha$, 
$hy_\alpha = g a_\alpha$ and $hx=x$ for each $x\in X\setminus \{ x_\alpha, y_\alpha, a_\alpha, g a_\alpha : \alpha<\lambda \}$. 
Since $h$ is macro-uniform, there exists $K\in \mathcal{E}$ such that $(x, y)\in H$ implies 
$(hx, hy)\in K$.
Hence, $(a_\alpha , g a_\alpha)\in K$ for each $\alpha< \lambda$.

If $| X \setminus (A_1\cup g A_1)| < |X|$ then we partition $A_1 = B\cup C$, $|B| = |C| = |X|$ and, to choose $K$, apply above arguments for the pair  $B, g B$ and  $C, g C$.

Repeating above construction for $A_2$ and $A_3$,
we see that $g$ is compatible with $\mathcal{E}$.

\vskip 7pt

$(ii)\Longrightarrow (i)$.
Let $G= \mathcal{S}_X ^{<\kappa}$,
$(X, \mathcal{E})= X_G$. We take an arbitrary $h\in \mathcal{S}_X$ and show that $h : X_G \rightarrow X_G $
is macro-uniform.

Let $F$  be a finite subset of 
$ \mathcal{S}_X ^{<\kappa}$,
$x\in X$, $y=hx$.
Then $hFh^{-1}\subset  \mathcal{S}_X ^{<\kappa}$ and,
for $f\in F$, we have 
$(hx, hfx)=(y, hfh^{-1}y)$, so $h$ is macro-uniform.
$ \ \ \Box$


\begin{thebibliography}{10}


\bibitem{b1} D. Dikranjan,  I.  Protasov, K. Protasova,  N. Zava,  {\em Balleans, hyperballeans and ideals}, Appl. Gen. Topology {\bf 20} (2019), 431-447.

\bibitem{b2} { I. V Protasov,} {\em   Balleans of bounded geometry  and  $G$-spaces},  Algebra Discrete Math.  {\bf 7}: 2 (2008), 101-108.

\bibitem{b3}{I. Protasov, } {\em Decompositions of set-valued mappings, } Algebra Discrete Math. {\bf 30:} 2 (2020), 235-238.

\bibitem{b4}{ I. Protasov, T. Banakh, }{\it Ball Structures and Colorings of Groups and Graphs},  Math. Stud. Monogr. Ser., vol. 11, VNTL, Lviv, 2003.


\bibitem{b5}{ I. Protasov, M. Zarichnyi,} {\em General Asymptology},   Math. Stud. Monogr. Ser., Vol. 12, VNTL, Lviv, 2007.


\bibitem{b6}    {J. Roe,} {\em Lectures on Coarse Geometry}, Univ. Lecture Ser., vol. 31, American Mathematical Society, Providence RI, 2003.


\end{thebibliography}

\vskip 15pt

CONTACT INFORMATION
\vskip 15pt

I.~Protasov: \\
Faculty of Computer Science and Cybernetics  \\
        Taras Shevchenko National University of Kyiv \\
         Academic Glushkov pr. 4d  \\
         03680 Kyiv, Ukraine \\ i.v.protasov@gmail.com

\end{document}